\documentclass[11pt]{article}
\usepackage{amsmath,amsthm,epsfig,amsfonts,color, bm}
\usepackage{rotating}
\allowdisplaybreaks

\usepackage{comment}
\usepackage{booktabs}
\usepackage{setspace}
\usepackage{array}

\newtheorem{theorem}{{Theorem}}
\newtheorem{lemma}{{Lemma}}

\usepackage{comment}

\newcommand{\E}{\mathbf{E}}
\renewcommand{\P}{\mathbf{P}}

\newcommand{\e}{\end{document}}

\begin{document}

\title{A Simplified Condition For Quantile Regression}

\author{Liang Peng\thanks{Maurice R. Greenberg School of Risk Science, Georgia State University, USA.} ~~~
and 
Yongcheng Qi
\thanks{Department of Mathematics and Statistics, University of Minnesota Duluth, USA.}
}

\date{}
\maketitle
\bigskip
\begin{quote}
\begin{center}
\textbf{Abstract}
\end{center}
Quantile regression is effective in modeling and inferring the conditional quantile given some predictors and has become popular in risk management due to wide applications of quantile-based risk measures. When forecasting risk for economic and financial variables, quantile regression has to account for heteroscedasticity, which raises the question of whether the identification condition on residuals in quantile regression is equivalent to one independent of heteroscedasticity. In this paper,  we present some identification conditions  under three probability models and use them to establish simplified conditions in quantile regression.

\smallskip


\end{quote}


\section{Introduction}
Since the distribution function and quantile function are two fundamental quantities in characterizing randomness, nonparametric distribution and quantile estimation have played a significant role in nonparametric statistics; see Shorlack and Wellner (1986) for an overview of empirical processes and quantile processes. Like  regression for modeling the relationship between response and predictors, quantile regression is an effective technique for modeling and inferring the conditional quantile of a response given predictors. We refer to Koenker (2005) for an overview of quantile regression techniques. Because quantile-based risk measures such as Value-at-Risk are popular in risk management, quantile regression plays an important role in forecasting risk. Below are three particular applications of using quantile regression to forecast risk in insurance and econometrics.

The first example is risk forecast in non-life insurance, where an actuarial dataset often includes the number of claims, the loss of each claim, and some characteristics,  such as the policyholder's age, type of car, and driving experience in automobile insurance, for each policyholder.  When a risk manager needs to forecast Value-at-Risk (VaR), a two-step inference procedure is employed in the literature: logistic regression for modeling the probability of having nonzero claims and quantile regression for modeling Value-at-Risk at an adjusted risk level computed from the logistic regression; see Heras,  Moreno and  Vilar-Zan\'{o}n (2018),  Kang, Peng and Golub (2021), and Fung et al. (2024). Because of many zero claims, the first step of logistic regression improves the inference for the probability of nonzero claims, while the second step of quantile regression effectively models and infers the quantile-based risk measure, leading to accurate risk forecasts. Also, the analysis often assumes independence across policyholders, which is practically sound.

The second example is about predictability in financial econometrics. Because of heterogeneous quantiles (see Gebka and Wohar (2013)  and  Ma, Xiao and Ma (2018)), 
researchers in econometrics are interested in using quantile regression to test for the predictability of some economic predictors; see
 Xiao (2009) for unit root predictor,  Xu (2021) for the case of stationary predictor, and
Lee (2016), Fan and Lee (2019), and Liu et al. (2023) for some unified tests regardless of predictors being stationary, near unit root, and unit root.

The third example is systemic risk, which has been a major concern in the financial industry and insurance business since the 2008 global financial crisis. A popular systemic risk measure is the so-called CoVaR defined as the conditional quantile of system loss at risk level q given some predictors and that an individual loss equals its Value-at-Risk (VaR) at the same risk level q (see Adrian and Brunnermeier (2016)).  Various applications and extensions of the systemic risk measure CoVaR have appeared in the literature as in Giglio, Kelly and Pruitt (2016), Yang and Hamori (2021), Girardi and Ergün (2013), H\"ardle, Wang and Yu (2016),  Chen, H\"ardle and Okhrin (2019), and Capponi and Rubtsov (2022).
Because CoVaR is defined as a conditional quantile given some economic predictors,
Adrian and Brunnermeier (2016) propose to model and infer by two quantile regression models: one quantile regression for individual loss given some macroeconomic predictors and another quantile regression for system loss given macroeconomic predictors and the individual loss.

A general quantile regression model for a univariate predictor is
\begin{equation}\label{QR1}
\left\{\begin{array}{ll}
Y_t=\alpha_q+{\beta}_q{X}_t+U_t,~ X_t=     {\mu+\rho X_{t-1}+e_t,~e_t=}\sum_{j=0}^{\infty}c_{j}V_{t-j},\\
\{(U_t, V_t)\}~\text{is a sequence of independent and identically distributed random vectors},
\end{array}\right.
\end{equation}
where $c_j$'s satisfy that $\{     {e_t}\}$ is stationary,      {$|\rho|<1$ for stationary $\{X_t\}$, and $\rho=1$ for the unit root process of $\{X_t\}$.}  Here, we allow the dependence between $U_t$ and $V_t$.
For $q\in (0, 1)$, to ensure that $\alpha_q+\beta_qX_t$ in (\ref{QR1}) models the q-th conditional quantile of $Y_t$ given $X_t$ and the consistency of quantile regression estimation, one commonly assumes that
\begin{equation}\label{con1}
\P(U_t\le0|X_t)=q.
\end{equation}
Then, a simple question is,  under model (\ref{QR1}), whether   (\ref{con1}) is equivalent to
\begin{equation}\label{con2}
\P(U_t\le0|V_t)=q.
\end{equation}

When quantile regression is applied to economic and financial variables, as in the second and third examples above, it is necessary to account for heteroscedasticity. In this case, one often considers the following quantile regression model:
\begin{equation}\label{QR2}
\left\{\begin{array}{ll}
Y_t=\alpha_q+{\beta}_q{X}_t+U_t,~ X_t=     {\mu+\rho X_{t-1}+e_t,~e_t=}\sum_{j=0}^{\infty}c_{j}V_{t-j},~V_t=\sigma_{t,x}\eta_t,~U_t=\sigma_{t,y}\varepsilon_t,\\
\{(\varepsilon_t, \eta_t)\}~\text{is a sequence of independent and identically distributed random vectors},\\
(\varepsilon_t, \eta_t)~ \text{is independent of}~ \{(\sigma_{s,y}, \sigma_{s,x}): s\le t\}, 
\end{array}\right.
\end{equation}
where $\{\sigma_{t,y}>0\}$ and $\{\sigma_{t,x}>0\}$ are stationary,  $c_j$'s satisfy that $\{     {e_t}\}$ is stationary,      {and $\rho\in (-1, 1]$. Again, $\rho\in (-1, 1)$ and $\rho=1$ represent stationary and nonstationary $\{X_t\}$, respectively}. Here, the second question is, under model (\ref{QR2}), whether   (\ref{con2}) is equivalent to
\begin{equation}\label{con3}
\P(\varepsilon_t\le 0|\eta_t)=q,
\end{equation}
     {which is equivalent to $\P(U_t\le 0|\eta_t)=q$ as $\sigma_{t,y}>0$, implying that the conditional quantile of $Y_t$ given $X_t$ is still $\alpha_q+\beta_qX_t$.}
Further,  the third question is whether  (\ref{con1}) is equivalent to (\ref{con3}) under model (\ref{QR2}). 

     {In using models (\ref{QR1}) and (\ref{QR2}) for the above mentioned three applications, $Y_t$ and $X_t$ are insurance loss and policyholder's characteristics respectively in forecasting insurance risk and
 are asset return and some economic variables respectively in predictability tests and systemic risk management.}
  
In this note, under some minor conditions, we prove the equivalence of (\ref{con1}) and (\ref{con2}) under model (\ref{QR1}) and the equivalence of (\ref{con2}) and (\ref{con3}) under model (\ref{QR2}). It is clear that (\ref{con3}) implies (\ref{con1}) under model (\ref{QR2}). Unfortunately, we can only show that  (\ref{con1}) implies (\ref{con3}) under model (\ref{QR2}) and some restrictive moment conditions.

\section{Main Results}

In this section, we will prove three technical lemmas that are of independent interest in probability theory. Then, we will prove three theorems that can be used to answer the above three questions.  

     {Throughout, we will use either characteristic function or moment generating function as they
determine distribution functions and work nicely with sums.}

\begin{lemma}\label{lem1} Assume random variables $U$ and $V$ are independent of $S$, and $S+U$
 and $S+V$ are identically distributed.  Then $U$ and $V$ have the same distribution under each of the following two  conditions:\\
 (a). The characteristic function of $S$, defined as  $\phi_S(t)=\E(e^{itS})$, is nonzero for $t\in C$, where $C$ is a dense subset of $(-\infty, \infty)$, where $i$ is the imaginary number with $i^2=-1$. \\ 
 (b). For some constant $h>0$, one of the two moment-generating functions
 $M_U(t):=\E\big(e^{tU}\big)$
 and $M_V(t):=\E\big(e^{tV}\big)$ is finite for $t\in (-h,h)$. 
\end{lemma}

\begin{proof} By using independence and equality in distributions we have, for all $t\in (-\infty, \infty)$,
\begin{equation}\label{eqofch}
\phi_S(t)\phi_U(t)=\E\big(e^{it(S+U)}\big)=
\E\big(e^{it(S+V)}\big)=\phi_S(t)\phi_V(t).
\end{equation}

Under condition (a) that $\phi_S(t)\ne 0$ for $t\in C$, we have from \eqref{eqofch} that 
$\phi_U(t)=\phi_V(t)$ for $t\in C$. Since all characteristic functions are continuous and $C$ is dense in $(-\infty, \infty)$, we have
$\phi_U(t)=\phi_V(t)$ for $t\in (-\infty, \infty)$, which implies that $U$ and $V$ are identically distributed.  

Under condition (b),  assume $M_U(t)$ is finite for $t\in (-h,h)$, implying that all moments of $U$ exist. It is known that when the moment-generating function of a random variable exists, the moments of the random variable uniquely determine the moment-generating function and the distribution of the random variable as well. In other words,  
if all moments of random variable $V$ exist and 
$\E(V^k)=\E(U^k)$ for all positive integers $k\ge 1$, then $U$ and $V$ are identically distributed. 

Now we will show $\E(V^k)=\E(U^k)$ for all $k\ge 1$.
Since $\phi_S(t)$ is continuous and $\phi_S(0)=1$, $\phi_S(t)\ne 0$ in a neighborhood of zero. Then it follow from \eqref{eqofch} that $\phi_U(t)=\phi_V(t)$, 
$t\in (-\delta, \delta)$ for some $\delta>0$. Hence,  both $\phi_U(t)$ and $\phi_V(t)$ are differentiable at the origin infinitely many times since all moments of $U$ exist, and  $\E(U^k)=(-i)^k\phi_U^{(k)}(t)|_{t=0}=(-i)^k\phi_V^{(k)}(t)|_{t=0}= \E(V^k)$ for all $k\ge 1$. This completes the proof.
\end{proof}

\begin{lemma}\label{lem2} Assume random variables $Y_0$ and $Y_1$ are independent of $Z$,  $\P(Z>0)=1$, and $ZY_0$
 and $ZY_1$ are identically distributed.  Then $Y_0$ and $Y_1$ have the same distribution under one of the following conditions.\\
(A). The characteristic function of $\ln Z$, $\phi_{\ln Z}(t)$, is nonzero for $t\in C$, where $C$ is a dense subset of $(-\infty, \infty)$. \\ 
(B). For some $\delta>0$,  $\E\big(|Y_j|^\delta+|Y_j|^{-\delta}I(Y_j\ne 0)\big)<\infty$ for $j=0$ or $1$.
\end{lemma}

\begin{proof}
Without loss of generality,  assume $p=\P(Y_0>0)>0$ and  
$r=\P(Y_0<0)>0$.
Then
$
\P(Y_1>0)=\P(ZY_1>0)=\P(ZY_0>0)=\P(Y_0>0)=p
$
since $\P(Z>0)=     {1}$.  Similarly,  we have
$\P(Y_0<0)=\P(Y_1<0)=r~ \text{and}~ \P(Y_0=0)=\P(Y_1=0)=1-p-r.$

Now define two new random variables, $U_0$ and $U_1$, with the following properties: 

\noindent\textbf{(P1)}: For $j=0,1$, the distribution of $U_j$ is the same as the conditional distribution of $Y_j$ given $Y_j>0$, that is,
\begin{equation}\label{cdf}
 \P(U_j\le u)=\P(Y_j\le u|Y_j>0)=\P(0<Y_j\le u)/p, ~~u>0. 
\end{equation}

 \noindent\textbf{(P2)}: $(U_0, U_1)$ and $Z$ are independent

Note that $\P(U_j>0)=1$ for $j=0,1$.
Since $ZY_0$ and $ZY_1$ are identically distributed, we have
\begin{equation}\label{equal-cond}
 \P(ZY_0\le x|ZY_0>0)=\P(ZY_1\le x|ZY_1>0), ~~~-\infty<x<\infty.  
\end{equation}
Again, since $\P(Z>0)=1$, we have for $j=0,1$
\begin{eqnarray*}
\P(ZU_{     {j}}\le x)&=&\P(ZY_j\le x|Y_j>0)=\frac{\P(ZY_j\le x, Y_j>0)}{\P(Y_j>0)}\\
&=&\frac{\P(ZY_j\le x, ZY_{     {j}}>0)}{\P(ZY_j>0)}=\P(ZY_j\le x|ZY_j>0).
\end{eqnarray*}
It follows from \eqref{equal-cond} that 
$ZU_0$ and $ZU_1$ are identically distributed, or equivalently, $\ln Z+\ln U_0$ and $\ln Z+\ln U_1$ have the same distribution.  Note that $\ln Z$ and $(\ln U_0, \ln U_1)$ are independent. 

To apply Lemma~\ref{lem1} with $S=\ln(Z)$, $U=\ln(U_0)$ and $V=
\ln(U_1)$, we need to verify conditions (a) and (b).  Clearly,   condition (A) in Lemma~\ref{lem2} implies condition (a) in Lemma~\ref{lem1}.  

Now assume condition (B) holds with $j=0$ or $1$. Then we have
$
\E(U_j^\delta)=\frac{1}{p}\E\big(Y_j^\delta I(Y_j>0)\big)\le \frac1p\E\big(|Y_j|^\delta\big)<\infty 
$
and 
$
\E(U_j^{-\delta})=\frac{1}{p}\E\big(Y_j^{-\delta} I(Y_j>0)\big)\le \frac1p\E\big(|Y_j|^{-\delta}I(Y_j\ne 0)\big)<\infty.
$
By using Lyapunov’s inequality, we have for any $t\in (0,\delta)$
\[
\big(\E(U_j^{t})\big)^{1/t}\le \big(\E(U_j^\delta)\big)^{1/\delta}<\infty~\mbox{ and }~
\big(\E(U_j^{-t})\big)^{1/t}\le \big(\E(U_j^{-\delta})\big)^{1/\delta}<\infty,
\]
that is,
$
\E\big(e^{\pm t\ln(U_j)}\big)<\infty,
$
or equivalently
$
\E\big(e^{t\ln(U_j)}\big)<\infty,~~t\in (-\delta,\delta),
$
which implies condition (b) in Lemma~\ref{lem1}.  Therefore, it follows 
from Lemma~\ref{lem1} that $\ln U_0$ and $\ln U_1$ are identically distributed,  that is, $U_0$ and $U_1$ are identically distributed.  From \eqref{cdf} we obtain 
\begin{equation}\label{result1}
\P(0<Y_0\le u)=\P(0<Y_1\le u),~~~u>0.    
\end{equation}

Replacing $Y_0$ and $Y_1$ with $-Y_0$ and $-Y_1$, respectively, and repeating the same procedure as above, we have
$
\P(0<-Y_0\le v)=\P(0<-Y_1\le v),~~~v>0,
$
which is equivalent to 
\begin{equation}\label{result2}
\P(v\le Y_0<0)=\P(v\le Y_1<0),~~~v<0.
\end{equation}

By combining \eqref{result1},  \eqref{result2}, and $\P(Y_0=0)=\P(Y_1=0)$,   we conclude that $Y_0$ and $Y_1$ are identically distributed.  This completes the proof of the lemma. 
\end{proof}

\begin{lemma}\label{lem3}
    Assume $Y_0$, $Y_1$, $Z$, and $W$ are random variables,  $(Y_0, Y_1)$ is independent of  $(Z, W)$,  and $\P(Z>0)=1$. 
Assume all moments of random variables $Z$ and $W$ are finite,  and the moment-generating functions of $Y_0$ and $Y_1$, $M_{Y_j}(t)=\E(e^{tY_j})$, exist in $t\in (-h, h)$ for some constant $h>0$. If $Y_0Z+W$ and $Y_1Z+W$ are identically distributed, then $Y_0$ and $Y_1$ are identically distributed. 
\end{lemma}

\begin{proof}
Since the moment-generating functions $\E(e^{tY_0})$ and $\E(e^{tY_1})$ are well defined in $t\in (-h, h)$, all moments of $Y_0$ and $Y_1$ are finite.  In this case, if $\E(Y_0^n)=\E(Y_1^n)$ for all $n\ge 1$, then $Y_0$ and $Y_1$ have the same moment-generating functions, consequently they have the identical distribution functions.  Hence, it suffices to  show  
 that 
 $\E (Y_0^n)=\E(Y_1^n)~\text{ for all }~n\ge 1.$

Note that all moments of $Z$ and $W$ exist,  $\mu_{j,k}:=\E(Z^jW^k)$ is finite for all $j, k\ge 0$, where all indices $     {j}, k$ are integer-values. We have $\mu_{j,0}=\E(Z^j)>0$ for all $j\ge 0$ since $\P(Z>0)=0$. 

For $n\ge 1$, we have
\[
g_n(t):=\E\big(tZ+W\big)^n=\E\Big(\sum^n_{m=0}\binom{n}{m}t^mZ^mW^{n-m}\Big)=\sum^n_{m=0}\mu_{m,n-m}\binom{n}{m}t^m.
\]
Since $Y_0Z+W$ and 
 $Y_1Z+W$ are identically distributed, 
$
\E(g_n(Y_j))=\E\{\E\big((Y_jZ+W)^n|Y_{     {j}}\big)\}=E\big((Y_jZ+W)^n\big)
$
is the same for $j=0,1$, that is,
\[
\sum^n_{m=0}\mu_{m,n-m}\binom{n}{m}\E(Y_0^m)=
     \sum^n_{m=0}\mu_{m,n-m}\binom{n}{m}\E(Y_1^m).
     \]
Rewrite the above equation as
\begin{equation}\label{n-th-moment} 
    \mu_{n,0}\big(\E(Y_0^n)-\E(Y_1^n)\big)=
    \sum^{n-1}_{m=0}\mu_{m,n-m}\binom{n}{m}\big(\E(Y_1^m)-\E(Y_0^m)\big).
\end{equation}

Notice that $\E(Y_0^0)=1=\E(Y_1^0)$ and $\mu_{n,0}=\E(Z^n)>0$ for all $n\ge 1$.  By setting $n=1$ in \eqref{n-th-moment}, we have
$\E(Y_0)=\E(Y_1)$. Now assume $\E(Y_0^m)=\E(Y_1^m)$ for all $m=0,1, \cdots, n-1$, then from \eqref{n-th-moment} we conclude $\E(Y_0^n)=\E(Y_1^n)$.  This implies that $\E(Y_0^n)=\E(Y_1^n)$ for $n\ge 1$ by induction.     
\end{proof}


\vspace{10pt}

\noindent\textbf{Remark 1.}  It is well known that the values of the characteristic function of a random variable in a neighborhood of the origin cannot uniquely determine the distribution of the random variable.  That is why we assume the characteristic function is non-zero over a dense subset of all real numbers; see condition (a) in Lemma~\ref{lem1} and condition (A) in Lemma~\ref{lem2},      {which require no moment conditions}.  For most of the commonly used distributions      {such as normal distributions,  t-distributions, and distributions in Remark 2 below},  these conditions are valid.       {When the characteristic function is zero in a set with a positive Lebesgue measure, we impose on the existence of moment generating function of another random variables to ensure the equality of  the distribution functions of the two variables; see condition (b) on Lemma~\ref{lem1} and condition (B) in Lemma~\ref{lem2}. Meanwhile,  the assumptions on the moment generating functions in Lemmas~\ref{lem1} - \ref{lem3} cannot be weakened to the finiteness of all moments of the random variables since moments cannot uniquely determine a distribution function in general;  see, e.g. Example 30.2, page 398 in  Billingsley (1995). 
}

\noindent\textbf{Remark 2.}
We notice that Smith (1962) showed that the characteristic function of a non-negative random variable cannot vanish identically in an interval. Immediately, we can draw the following two conclusions:
\noindent\textbf{i)} The characteristic function of a random variable $S$ is non-zero over a dense subset of all real numbers if the random variable $S$ is bounded from above or bounded from below.
\noindent\textbf{ii)} For a positive random variable $Z$, if $\P(Z>c)=1$ for some $c>0$, then the characteristic function of $\ln Z$ is non-zero over a dense subset of all real numbers.

\noindent\textbf{Remark 3.}  In Lemma~\ref{lem3}, we impose moment conditions      {via moment generating functions as it remains unknown whether or how a condition on the characteristic function like Lemmas~\ref{lem1} and \ref{lem2} can be employed}.

\bigskip

First, we show Theorem \ref{thm1} below. The equivalence between (\ref{con1}) and (\ref{con2}) under model (\ref{QR1}) follows immediately.

\begin{theorem}\label{thm1}
Assume random variable $W$ is independent of random vector $(X, Y)$.  Assume constant $q\in (0, 1)$.  Then,  
$
\P(X\le0|Y+W)=q~\text{if and only if}~ \P(X\le0|Y)=q
$
under one of the following two conditions:

\noindent (C1):  The characteristic function of $W$, $\phi_W(t)=\E\big(e^{itW}\big)$, is non-zero in a dense subset of $(-\infty, \infty)$.   

\noindent (C2): The moment-generating function of $Y$, $M_Y(t)$,   is finite in $(-h,h)$ for some $h>0$. 
\end{theorem}

\begin{proof}
When $\P(X\le0|Y)=q$, we have
\begin{equation}\label{pfTh1-0}
\begin{array}{ll}
&\P(X\le0|Y+W)=\E\{\P(X\le0|Y, W)|Y+W\}\\
=&\E\{\P(X\le0|Y)|Y+W\}=\E(q|Y+W)=q.
\end{array}\end{equation}
Hence, we only need to prove $\P(X\le0|Y)=q$ when
\begin{equation}\label{pfTh1-1}
\P(X\le0|Y+W)=q.
\end{equation}
It follows from \eqref{pfTh1-1} that
$
\P(X\le0|Y+W)=q=\P(X\le 0)$ and $\P(X>0|Y+W)=1-q=\P(X>0),
$
which implies that the Bernoulli random variable $T:=I(X\le 0)$ is independent of $Y+W$. 

For $j=0,1$, define the conditional distribution of $Y$ given $T=j$ as
\begin{equation}\label{cond-cdf}
F_j(y)=\P(Y\le y|T=j)~~~\text{for } -\infty<y<\infty.
\end{equation}

Our objective is to show that $T$ and $Y$ are independent, or equivalently, the conditional distribution of $Y$ given $T=j$ is the same for $j=0,1$.  Now assume random variables $Y_0$ and $Y_1$ are independent of $W$, and their cumulative distribution functions are $F_0$ and $F_1$, respectively, as defined in \eqref{cond-cdf}.  For every $x$, we have
$
\P(W+Y_j\le x)=\P(W+Y\le x|T=j)~~~\text{for}~j=0, 1.
$
Since $T$ and $W+Y$ are independent, we have
$
\P(W+Y\le x|T=j)=\P(W+Y\le x)~~~\text{for}~j=0, 1.
$
Therefore, $W+Y_0$ and $W+Y_1$ are identically distributed. By applying Lemma~\ref{lem1} with $S=W$, $U=Y_0$ and $V=Y_1$, we conclude that $Y_0$ and $Y_1$ have the same distribution function if we can verify condition (a) or (b) in Lemma~\ref{lem1} holds.  In fact,  condition (C1) implies condition (a) in Lemma~\ref{lem1}. 
Under condition (C2), we have
$
M_{Y_1}(t)=\E\big(e^{tY_1}\big)=\frac{1}{q}
\E\big(e^{tY}I(X\le 0)\big)\le \frac{1}{q}M_{Y}(t)<\infty
$
for all $t\in (-h,h)$, that is, condition (b) holds in Lemma~\ref{lem1}.

Thus,  we have proved that $T$ and $Y$ are independent, which yields
$
\P(X\le 0|Y)=\P(T=1|Y)=\P(T=1)=q.
$
This completes the proof.
\end{proof}

Second, we prove Theorem \ref{thm2} below. The equivalence between (\ref{con2}) and (\ref{con3}) under model (\ref{QR2}) follows immediately by noting that $I(U_t\le 0)=I(\varepsilon_t\le 0)$.

\begin{theorem}\label{thm2}
Assume random variable $Z$ is independent of random vector $(X, Y)$, $\P(Z>0)=1$, and    $q\in (0, 1)$ is a constant. Then,   
\[
\P(X\le0|YZ)=q~\text{if and only if}~ \P(X\le0|Y)=q
\]
under one of the following two conditions:

\noindent (C3):  The characteristic function of $\ln Z$, $\phi_{\ln Z}(t)=\E\big(e^{it\ln Z}\big)$, is non-zero in a dense subset of $(-\infty, \infty)$.   

\noindent (C4):  For some $\delta>0$,  $\E\big(|Y|^\delta+|Y|^{-\delta}I(Y\ne 0)\big)<\infty$.
\end{theorem}

\begin{proof} The sufficiency follows from the same line as in the proof of (\ref{pfTh1-0}).  Therefore,  we only need to prove $\P(X\le 0|Y)=q$ when
\begin{equation}\label{pfTh2-1}
\P(X\le 0|YZ)=q.
\end{equation}

As in the proof of Theorem~\ref{thm1},  set $T=I(X\le 0)$ and define random variables $Y_0$ and $Y_1$ in such a way that $Y_0$ and  $Y_1$ are independent of $Z$ and corresponding distribution functions are $F_0$ and $F_1$ as defines in \eqref{cond-cdf}.   
It suffices to show the independence of $T$ and $Y$, which is equivalent to the equality of the distribution functions of $Y_0$ and $Y_1$.

Since \eqref{pfTh2-1} implies the independence of 
$T$ and $YZ$, we have for $x\in (-\infty, \infty)$
$
\P(ZY\le x|T=0)=\P(ZY\le x|T=1).
$
Because the left-hand side and the right-hand side of the above equation are
equal to $\P(ZY_0\le x)$ and $\P(ZY_1\le x)$, respectively,  we conclude that $ZY_0$ and $ZY_1$ are identically distributed. Clearly, condition (C3) is the same as condition (A) in Lemma~\ref{lem2}.  We can also show that condition  (C4) implies condition (B) in Lemma~\ref{lem2} since
\begin{eqnarray*}
\E\big(|Y_1|^\delta+|Y_1|^{-\delta}I(Y_1\ne 0)\big)
&=&
\frac{1}{q}\E\big\{\big(|Y|^\delta+|Y|^{-\delta}I(Y\ne 0)\big)I(T=0)\big\}\\
&\le &
\frac1q\E\big(|Y|^\delta+|Y|^{-\delta}I(Y\ne 0)\big).
\end{eqnarray*}
Hence, it follows from Lemma~\ref{lem2} that
$Y_0$ and $Y_1$ are identically distributed, that is, 
$\P(Y\le y|T=0)=\P(Y\le y|T=1)=\P(Y\le y)$,   
and $Y$ is independent of $I(X\le 0)$. Therefore,
$\P(X\le 0|Y)=\P(X\le0)=q.$
That is, Theorem~\ref{thm2} holds.
\end{proof}

Third, like the proof of (\ref{pfTh1-0}), we know that (\ref{con3}) implies (\ref{con1}) under model (\ref{QR2}). But, we can only show the equivalence between (\ref{con1}) and (\ref{con3}) under model (\ref{QR2}) with enough finite moments by using the following theorem.


\begin{theorem}\label{thm3}
Assume $X$, $Y$, $Z$, and $W$ are four random variables,  $(X, Y)$ and $(Z, W)$ are independent,  and $\P(Z>0)=1$. 
Furthermore, we assume all moments of random variables $Z$ and $W$ are finite,  and the moment-generating function of $Y$, $M_Y(t)=\E(e^{tY})$, exists in $t\in (-h, h)$ for some constant $h>0$.  Assume constant $q\in (0, 1)$. Then,  

\begin{equation}\label{C2}
\P(X\le 0|YZ+W)=q~\text{if and only if}~ \P(X\le 0|Y)=q.
\end{equation}
\end{theorem}

\begin{proof}  As before, if $\P(X\le 0|Y)=q$, then we have
\begin{eqnarray*}
\P(X\le 0|YZ+W)&=&\E\{\P(X\le 0|Y, Z, W)|YZ+W\}\\
&=&\E\{\P(X\le 0|Y)|YZ+W\}\\
&=&\E(q|YZ+W)=q.
\end{eqnarray*}

Now we assume  $\P(X\le 0|YZ+W)=q$. Follow the proof of Theorem~\ref{thm1}, we define random variable $T=I(X\le 0)$. Then
$
\P(T=1|YZ+W)=q$ and $\P(T=0|YZ+W)=1-q,
$
that is, $T$ and $YZ+W$ are independent and $\P(T=1)=q$.

Define distributions $F_0$ and $F_1$, as in \eqref{cond-cdf}. If $F_0$ and $F_1$ are the same, then we can conclude the independence of $T$ and $Y$, which implies
$\P(X\le 0|Y)=\P(T=1|Y)=\P(T=1)=q$.
Hence, we only need to show that $F_0$ and $F_1$ are identical. 

Like the proof of Theorem~\ref{thm1}, we assume random variables $Y_0$ and $Y_1$ are independent of $(Z, W)$, and their cumulative distribution functions are $F_0$ and $F_1$, respectively. Using the independence of $(X, Y)$ and $(Z,W)$, we have
for each $j=0,1$,
the distribution of $Y_jZ+W$ is the same as the conditional distribution of $YZ+W$ given $T=j$.  Meanwhile, the conditional distribution of $YZ+W$ given $T=j$ is the same for $j=0, 1$ due to the independence of   
$T$ and $YZ+W$.
Therefore, we know that $Y_0Z+W$ and 
 $Y_1Z+W$ are identically distributed.   
 Since
\[
\E(e^{tY})=\E(e^{tY}|T=0)(1-q)+\E(e^{tY}|T=1)q=(1-q)\E(e^{tY_0})+q\E(e^{tY_1})
\]
 for all $t\in (-h,h)$,  the moment-generating functions for $Y_0$ and $Y_1$ are well defined in $(-h, h)$. 
 Hence, an application of Lemma~\ref{lem3} concludes that $Y_0$ and $Y_1$ are identically distributed, that is,  $F_0=F_1$, i.e.,  the theorem holds.  
\end{proof}



\end{document}